# GOURSAT PROBLEM IN THE NON-CLASSICAL TREATMENT FOR A SIXTH ORDER PSEUDOPARABOLIC EQUATION

**I.G.Mamedov**

*A.I.Huseynov Institute of Cybernetics of NAS of Azerbaijan. Az 1141,
Azerbaijan, Baku st. B. Vahabzade, 9
E-mail: ilgar-mammadov@rambler.ru*

## Abstract

*In the paper the Goursat problem with non classical boundary conditions not requiring the agreement conditions is considered for a sixth order pseudoparabolic equation with the classical boundary condition is substantiated in the case if the solution of the stated problem is sought in S.L.Sobolev anisotropic space.*

## Problem statement

Consider equation

$$(V_{2,4}u)(x) \equiv D_1^2 D_2^4 u(x) + \sum_{\substack{i=0 \\ i+j<6}}^{2} \sum_{j=0}^{4} a_{i,j}(x) D_1^i D_2^j u(x) = Z_{2,4}(x) \in L_p(G) \qquad (1)$$

Here $u(x) \equiv u(x_1, x_2)$ is a desired function determined on $G$; $a_{i,j}(x)$ are the given measurable functions on $G = G_1 \times G_2$, where $G_k = (0, h_k)$, $k = \overline{1,2}$; $Z_{2,4}(x)$ is a given measurable function on $G$; $D_k^\xi = \partial^\xi / \partial x_k^\xi$ is a generalized differentiation operator in S.L.Sobolev sense, $D_k^0$ is an identity transformation operator.

Equation (1) is a hyperbolic equation and has two characteristics $x_1 = const, x_2 = const$, one of which is double-fold the second one is four-fold. Therefore, in some sense we can consider equation (1) as a pseudoparabolic equation [1]. This equation is a Boussenesq-Liav generalized equation from the vibrations theory [2] and Aller's equation under mathematical modeling [3, p.261] of the moisture absorption process in biology.

In the present paper equation (1) is considered in the general case when the coefficients $a_{i,j}(x)$ are non-smooth functions satisfying only the following conditions:

$$a_{i,j}(x) \in L_p(G), \; i = \overline{0,1} \; j = \overline{0,3};$$

$$a_{2,j}(x) \in L_{\infty,p}^{x_1,x_2}(G), \; j = \overline{0,3};$$

$$a_{i,4}(x) \in L_{p,\infty}^{x_1,x_2}(G), \; i = \overline{0,1}.$$

Therewith, the important principal moment is that the considered equation possesses discontinuous coefficients satisfying only some $p$-integrability and boundedness conditions i.e. the considered pseudoparabolic operator $V_{2,4}$ has no traditional conjugated operator. In other words, the Riemann function for this equation can't be investigated by the classical method of characteristics. In the papers [4-5] The Riemann function is determined as the solution of an integral equation. This is more natural than the classical way for deriving the Riemann function. The matter is that in the classic variant, for determining the Riemann function, the rigid smooth conditions on the coefficients of the equation are required.

Under these conditions, we'll look for the solution $u(x)$ of equation (1) in S.L.Sobolev anisotropic space

$$W_p^{(2,4)}(G) \equiv \left\{ u(x) : D_1^i D_2^j u(x) \in L_p(G), \; i = \overline{0,2}, \; j = \overline{0,4} \right\},$$

where $1 \leq p \leq \infty$. We'll define the norm in the space $W_p^{(2,4)}(G)$ by the equality

$$\|u\|_{W_p^{(2,4)}(G)} \equiv \sum_{i=0}^{2} \sum_{j=0}^{4} \left\| D_1^i D_2^j u \right\|_{L_p(G)}$$

For equation (1) we can give the classic form Goursat condition in the form (see. Fig. 1):

$$\begin{cases} u(0,x_2) = \varphi_1(x_2); \quad u(x_1,0) = \psi_1(x_1); \\ \left.\dfrac{\partial u(x)}{\partial x_1}\right|_{x_1=0} = \varphi_2(x_2); \quad \left.\dfrac{\partial u(x)}{\partial x_2}\right|_{x_2=0} = \psi_2(x_1); \\ \left.\dfrac{\partial^2 u(x)}{\partial x_2^2}\right|_{x_2=0} = \psi_3(x_1); \quad \left.\dfrac{\partial^3 u(x)}{\partial x_2^3}\right|_{x_2=0} = \psi_4(x_1), \end{cases} \qquad (2)$$

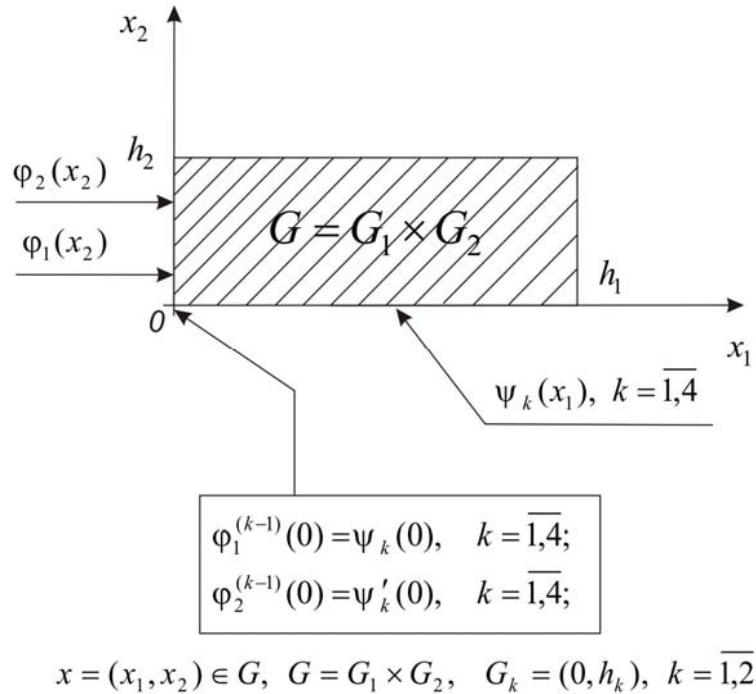

$$\varphi_1^{(k-1)}(0) = \psi_k(0), \quad k = \overline{1,4};$$
$$\varphi_2^{(k-1)}(0) = \psi_k'(0), \quad k = \overline{1,4};$$

$$x = (x_1, x_2) \in G, \quad G = G_1 \times G_2, \quad G_k = (0, h_k), \quad k = \overline{1,2}$$

**Fig.1. Geometric interpretation of Goursat classical boundary conditions.**

where $\varphi_1(x_2)$, $\varphi_2(x_2)$ and $\psi_k(x_1)$, $k = \overline{1,4}$ are the given measurable functions on $G$. It is obvious that in the case of conditions (2), in addition to the conditions

$$\varphi_1(x_2) \in W_p^{(4)}(G_2), \quad \varphi_2(x_2) \in W_p^{(4)}(G_2)$$

and

$$\psi_k(x_1) \in W_p^{(2)}(G_1), \quad k = \overline{1,4}$$

the given functions should also satisfy the following agreement conditions:

$$\begin{cases} \varphi_1(0) = \psi_1(0); & \varphi_2(0) = \psi_1'(0) \\ \varphi_1'(0) = \psi_2(0); & \varphi_2'(0) = \psi_2'(0) \\ \varphi_1''(0) = \psi_3(0); & \varphi_2''(0) = \psi_3'(0) \\ \varphi_1'''(0) = \psi_4(0); & \varphi_2'''(0) = \psi_4'(0). \end{cases} \qquad (3)$$

Consider the following non-classical initial-boundary conditions:

$$\begin{cases} V_{i,j}u \equiv D_1^i D_2^j u(0) = Z_{i,j} \in R, i = \overline{0,1}, j = \overline{0,3}; \\ (V_{2,j}u)(x_1) \equiv D_1^2 D_2^j u(x_1,0) = Z_{2,j}(x_1) \in L_p(G_1), j = \overline{0,3}; \\ (V_{i,4}u)(x_2) \equiv D_1^i D_2^4 u(0, x_2) = Z_{i,4}(x_2) \in L_p(G_2), i = \overline{0,1}; \end{cases} \qquad (4)$$

If the function $u \in W_p^{(2,4)}(G)$ is a solution of the classical form Goursat problem (1), (2), then it is also a solution of problem (1), (4) for $Z_{i,j}$, defined by the following equalities:

$$Z_{0,0} = \psi_1(0) = \varphi_1(0); \qquad Z_{1,0} = \varphi_2(0) = \psi_1'(0);$$

$$Z_{0,1} = \psi_2(0) = \varphi_1'(0); \qquad Z_{1,1} = \varphi_2'(0) = \psi_2'(0);$$

$$Z_{0,2} = \psi_3(0) = \varphi_1''(0); \qquad Z_{1,2} = \varphi_2''(0) = \psi_3'(0);$$

$$Z_{0,3} = \psi_4(0) = \varphi_1'''(0) \qquad Z_{1,3} = \varphi_2'''(0) = \psi_4'(0);$$

$$Z_{2,0}(x_1) = \psi_1''(x_1); \qquad Z_{2,1}(x_1) = \psi_2''(x_1)$$

$$Z_{2,2}(x_1) = \psi_3''(x_1); \qquad Z_{2,3}(x_1) = \psi_4''(x_1);$$

$$Z_{0,4}(x_2) = \varphi_1^{(IV)}(x_2); \qquad Z_{1,4}(x_2) = \varphi_2^{(IV)}(x_2)$$

The inverse one is easily proved. In other words, if the function $u \in W_p^{(2,4)}(G)$ is a solution of problem (1), (4), then it is also a solution of problem (1), (2) for the following functions:

$$\varphi_1(x_2) = Z_{0,0} + x_2 Z_{0,1} + \frac{x_2^2}{2!} Z_{0,2} + \frac{x_2^3}{3!} Z_{0,3} + \int_0^{x_2} \frac{(x_2 - \tau)^3}{3!} Z_{0,4}(\tau) d\tau; \qquad (5)$$

$$\varphi_2(x_2) = Z_{1,0} + x_2 Z_{1,1} + \frac{x_2^2}{2!} Z_{1,2} + \frac{x_2^3}{3!} Z_{1,3} + \int_0^{x_2} \frac{(x_2 - \xi)^3}{3!} Z_{1,4}(\xi) d\xi; \qquad (6)$$

$$\psi_1(x_1) = Z_{0,0} + x_1 Z_{1,0} + \int_0^{x_1}(x_1 - \eta)Z_{2,0}(\eta)d\eta; \qquad (7)$$

$$\psi_2(x_1) = Z_{0,1} + x_1 Z_{1,1} + \int_0^{x_1}(x_1 - \nu)Z_{2,1}(\nu)d\nu; \qquad (8)$$

$$\psi_3(x_1) = Z_{0,2} + x_1 Z_{1,2} + \int_0^{x_1}(x_1 - \mu)Z_{2,2}(\mu)d\mu; \qquad (9)$$

$$\psi_4(x_1) = Z_{0,3} + x_1 Z_{1,3} + \int_0^{x_1}(x_1 - \tau)Z_{2,3}(\tau)d\tau; \qquad (10)$$

Note that the functions (5)-(10) possess one important property, more exactly, for all $Z_{i,j}$, the agreement conditions (3) possessing the above-mentioned properties are fulfilled for them automatically. Therefore, equalities (5)-(10) may be considered as a general kind of all the functions $\varphi_1(x_2), \varphi_2(x_2)$ and $\psi_k(x_1)$ $k = \overline{1,4}$ satisfying the agreement conditions (3).

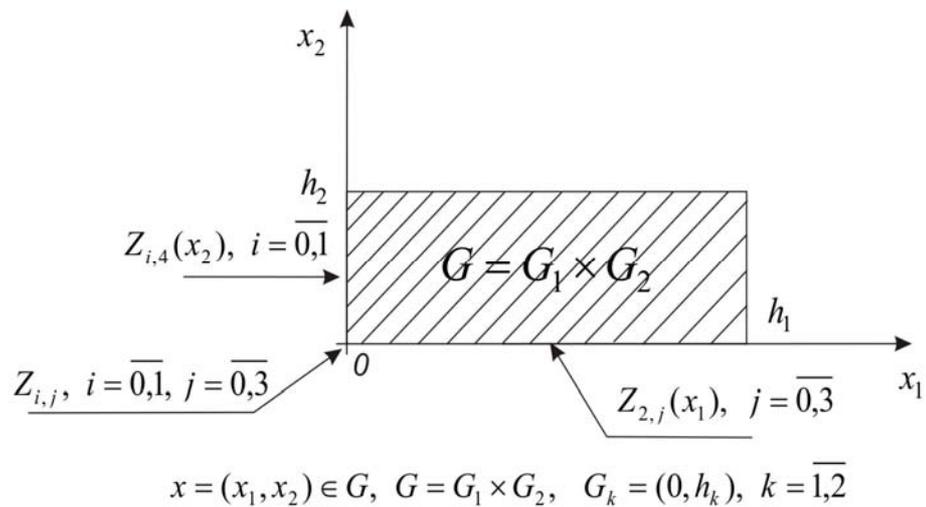

Fig.2. Geometrical interpretation of Goursat
boundary conditions in non – classical statement

So, the classical form Goursat problems (1), (2) and in non classical treatment (1), (4) (see. fig. 2) are equivalent in the general case. However, the Goursat problem in non-classical statement (1), (4) is more natural by statement than problem (1), (2). This is connected with the fact that in statement of problem (1), (4) the right sides of boundary conditions don't require additional conditions of agreement type. Note that some Goursat problems in non-classical treatments for hyperbolic and also pseudoparabolic equations were investigated in the author's papers [6-9].